\documentclass[leqno, 12pt]{amsart}
\usepackage{amscd}
\usepackage{amssymb}

\newcommand{\nc}{\newcommand}
\newcommand{\rc}{\renewcommand}

\nc{\ra}{{      \rightarrow     }}
\nc{\laa}{{     \leftarrow      }}      
\nc{\lra}{{\longrightarrow}}
\nc{\lr}{{\leftrightarrow}}             
\nc{\lrs}{{\rightleftarrows}}           

\nc{\imp}{{\Rightarrow}}                
\nc{\eq}{{\Leftrightarrow}}             

\nc{\inj}{{\pr  \hookrightarrow }}              
\nc{\injj}{{\pr \hookleftarrow  }}              

\nc{\sur}{{     \twoheadrightarrow      }}      
\nc{\surr}{{    \twoheadleftarrow       }}      
                        
\nc{\mm}{{\mapsto}}                             

\nc{\va}{{\uparrow}}                            

\nc{\barr}{\overline}
\nc{\ul}{\underline}
\nc{\sub}{\subseteq}

\nc{\se}{       \section                }
\nc{\sus}{      \subsection             }
\nc{\sss}{      \subsubsection          }

\nc{\Lemm}{     \subsection{Lemma}              }
\nc{\lemm}{     \subsubsection{Lemma}           }
\nc{\slemm}{    \subsubsection*{Lemma}          }
\nc{\sublemm}{  \subsubsection{\bf Sublemma}            }
\nc{\ssublemm}{         \subsubsection*{\bf Sublemma}           }

\nc{\Pro}{      \subsection{Proposition}        }
\nc{\pro}{      \subsubsection{Proposition}     }
\nc{\spro}{     \subsubsection*{Proposition}    }

\nc{\Corr}{     \subsection{Corollary}          }
\nc{\corr}{     \subsubsection{Corollary}       }
\nc{\scorr}{    \subsubsection*{Corollary}      }

\nc{\Theo}{     \subsection{Theorem}            }               
\nc{\theo}{ 	\subsubsection{\em Theorem}		}
\nc{\stheo}{    \subsubsection*{Theorem}        }

\nc{\rem}{      \subsubsection{Remark}          }
\nc{\srem}{     \subsubsection*{Remark} }

\nc{\rems}{     \subsubsection{Remarks}         }
\nc{\srems}{    \subsubsection*{Remarks}        }

\nc{\conj}{     \subsubsection{Conjecture}      }
\nc{\sconj}{    \subsubsection*{Conjecture}     }

\nc{\ex}{       \subsubsection{Example}         }
\nc{\sex}{      \subsubsection*{Example}        }
\nc{\exs}{      \subsubsection{Examples}        }
\nc{\sexs}{     \subsubsection*{Examples}       }

\nc{\que}{      \subsubsection{Question}        }
\nc{\ques}{     \subsubsection{Questions}       }
\nc{\sque}{     \subsubsection*{Question}       }
\nc{\sques}{    \subsubsection*{Questions}      }

\nc{\pl}{{\oplus}}                              
\nc{\tim}{{\times}}             
\nc{\btim}{{\boxtimes}}
\nc{\ltim}{\ltimes}                     %
\nc{\rtim}{\rtimes}                     %
\nc{\ltr}{\triangleleft}        %
\nc{\rtr}{\triangleright}       %

\nc{\ten}{{     \otimes         }}            
\nc{\Lten}{{    \aa{L}\otimes   }}            
\nc{\Ltim}{{    \aa{L}\times    }}            
\nc{\bten}{{\boxtimes}}                         

\nc{\con}{{ @>{\protect\cong}>> }}      
\nc{\conl}{{    @>{\cong}>>     }}      
\nc{\conn}{{    @<{\cong}<<     }}      
\nc{\Con}{{     \equiv          }}      
\nc{\appr}{{    \sim            }}      
\nc{\eqr}{{     \sim            }}      

\nc{\ha}{{ \frac{1}{2} }}               
        \nc{\half}{{ \frac{1}{2} }}

\nc{\ci}{{\circ}}               
\nc{\cd }{{\cdot}}              
\nc{\cddd}{{\cdots}}

\nc{\cupp}{\bigcup}             
\nc{\capp}{\bigcap}
\nc{\pll}{\bigoplus}

\nc{\pii}{\prod}                
\nc{\ppii}{\bigprod}            
\nc{\cci}{\sqcup}              
\nc{\ccii}{\bigsqcup}

\nc{\wwe}{\bigwedge}            
\nc{\cce}{\bigcoprod}           

\nc{\aaa}{      \stackerel      }       

\rc{\Im}{{      \operatorname{Im}       }}
\nc{\rank}{{    \ \operatorname{rank}\  }}
\nc{\Res}{{     \  \operatorname{Res}   }}
\nc{\Hom}{{    \operatorname{Hom}      }}
\nc{\End}{{     \operatorname{End}      }}
\nc{\RHom}{{    \operatorname{RHom}     }}
\nc{\HHom}{{    \operatorname{$\HH$om}  }}
\nc{\EEnd}{{    \operatorname{$\EE nd$} }}
\nc{\AAut}{{    \operatorname{$\AA ut$} }}
\nc{\RHHom}{{   \operatorname{R$\HH$om} }}
\nc{\Ext}{\operatorname{Ext}}
\nc{\Der}{{     \operatorname{Der}      }}

\nc{\ord        }{{ \operatorname{ord} }}                       
\nc{\divv       }{{ \operatorname{div} }}                       
\nc{\Lie        }{{ \operatorname{Lie} }}

\rc{\AA}{{\mathcal A}}
\nc{\BB}{{\mathcal B}} 
\nc{\CC}{{\mathcal C}}
\nc{\DD}{{\mathcal D}}
\nc{\EE}{{\mathcal E}}
\nc{\FF}{{\mathcal F}}
\nc{\GG}{{\mathcal G}}
\nc{\HH}{{\mathcal H}}
\nc{\II}{{\mathcal I}}
\nc{\JJ}{{\mathcal J}}
\nc{\KK}{{\mathcal K}}
\nc{\LL}{{\mathcal L}}
\nc{\MM}{{\mathcal M}}
\nc{\NN}{{\mathcal N}}
\nc{\OO}{{\mathcal O}}
\nc{\PP}{{\mathcal P}}
\nc{\QQ}{{\mathcal Q}}
\nc{\RR}{{\mathcal R}}
\rc{\SS}{{\mathcal S}}
\nc{\UU}{{\mathcal U}}
\nc{\VV}{{\mathcal V}}
\nc{\WW}{{\mathcal W}}
\nc{\ZZ}{{\mathcal Z}}
\nc{\XX}{{\mathcal X}}
\nc{\YY}{{\mathcal Y}}

\nc{\A}{{\mathbb A }}
\nc{\B}{{\mathbb B}}
\nc{\C}{{\mathbb C}}
\nc{\D}{{\mathbb D}}
\nc{\E}{{\mathbb E}}
\nc{\F}{{\mathbb F}}
\nc{\G}{{\mathbb G}}
\nc{\hH}{{\mathbb H}}
\nc{\I}{{\mathbb I}}
\nc{\J}{{\mathbb J}}
\nc{\K}{{\mathbb K}}
\nc{\lL}{{\mathbb L}}
\nc{\M}{{\mathbb M}}
\nc{\N}{{\mathbb N}}
\nc{\oO}{{\mathbb O}}
\nc{\pP}{{\mathbb P}}      
\nc{\Q}{{\mathbb Q}}
\nc{\R}{{\mathbb R}}
\nc{\sS}{{\mathbb S}}
\nc{\T}{{\mathbb T}}
\nc{\U}{{\mathbb U}}
\nc{\V}{{\mathbb V}}
\nc{\W}{{\mathbb W}}
\nc{\Z}{{\mathbb Z}}
\nc{\X}{{\mathbb X}}
\nc{\Y}{{\mathbb Y}}

\nc{\fA}{{\mathfrak A}}
\nc{\fB}{{\mathfrak B}}
\nc{\fC}{{\mathfrak C}}
\nc{\fD}{{\mathfrak D}}
\nc{\fE}{{\mathfrak E}}
\nc{\fF}{{\mathfrak F}}
\nc{\fG}{{\mathfrak G}}
\nc{\fH}{{\mathfrak H}}
\nc{\fI}{{\mathfrak I}}
\nc{\fJ}{{\mathfrak J}}
\nc{\fK}{{\mathfrak K}}
\nc{\fL}{{\mathfrak L}}
\nc{\fM}{{\mathfrak M}}
\nc{\fN}{{\mathfrak N}}
\nc{\fO}{{\mathfrak O}}
\nc{\fP}{{\mathfrak P}}
\nc{\fQ}{{\mathfrak Q}}
\nc{\fR}{{\mathfrak R}}
\nc{\fS}{{\mathfrak S}}
\nc{\fT}{{\mathfrak T}}
\nc{\fU}{{\mathfrak U}}
\nc{\fV}{{\mathfrak V}}
\nc{\fW}{{\mathfrak W}}
\nc{\fZ}{{\mathfrak Z}}
\nc{\fX}{{\mathfrak X}}
\nc{\fY}{{\mathfrak Y}}
\nc{\fa}{{\mathfrak a}}
\nc{\fb}{{\mathfrak b}}
\nc{\fc}{{\mathfrak c}}
\nc{\fd}{{\mathfrak d}}
\nc{\fe}{{\mathfrak e}}
\nc{\ff}{{\mathfrak f}}
\nc{\fg}{{\mathfrak g}}
\nc{\fh}{{\mathfrak h}}
\nc{\fj}{{\mathfrak j}}
\nc{\fk}{{\mathfrak k}}
\nc{\fl}{{\mathfrak{l}}}
\nc{\fm}{{\mathfrak m}}
\nc{\fn}{{\mathfrak n}}
\nc{\fo}{{\mathfrak o}}
\nc{\fp}{{\mathfrak p}}
\nc{\fq}{{\mathfrak q}}
\nc{\fr}{{\mathfrak r}}
\nc{\fs}{{\mathfrak s}}
\nc{\ft}{{\mathfrak t}}
\nc{\fu}{{\mathfrak u}}
\nc{\fv}{{\mathfrak v}}
\nc{\fw}{{\mathfrak w}}
\nc{\fz}{{\mathfrak z}}
\nc{\fx}{{\mathfrak x}}
\nc{\fy}{{\mathfrak y}}

\nc{\al}{{\alpha }}
\nc{\be}{{\beta }}
\nc{\ga}{{\gamma }}
\nc{\de}{{\delta }}
\nc{\del}{{\partial }}
\nc{\ep}{{\varepsilon }}
\nc{\vap}{{\epsilon }}

\nc{\ze}{{\zeta }}
\nc{\et}{{\eta }}
\rc{\th}{{\theta }}
\nc{\vth}{{\vartheta }}

\nc{\io}{{\iota }}
\nc{\ka}{{\kappa }}
\nc{\la}{{\lambda }}
\nc{\vrho}{{\varrho}}
\nc{\si}{{\sigma }}
\nc{\ups}{{\upsilon }}
\nc{\vphi}{{\varphi }}
\nc{\om}{{\omega }}

\nc{\Ga}{{\Gamma }}
\nc{\De}{{\Delta }}
\nc{\nab}{{\nabla}}
\nc{\Th}{{\Theta }}
\nc{\La}{{\Lambda }}
\nc{\Si}{{\Sigma }}
\nc{\Ups}{{\Upsilon }}
\nc{\Om}{{\Omega }}

\nc{\toc}{{\tableofcontents}}
\nc{\addl}{     \addcontentsline{toc}{subsection}       }

\begin{document}

\nc{\GK}{{      G(\KK)          }}
\nc{\GO}{{      G(\OO)          }}

\nc{\Gd}{{  {\check G}          }}

\nc{\Gh}{{  \hat{\GG}           }}
\nc{\GA}{{  G(\AA)              }}
\nc{\BK}{{  B(\mathcal K)           }}
\nc{\BO}{{  B(\mathcal O)           }}
\nc{\BKz}{{  B(\mathcal K)_0        }}
\nc{\NK}{{  N(\mathcal K)           }}
\nc{\NO}{{  N(\mathcal O)           }}
\nc{\TK}{{  B(\mathcal K)           }}
\nc{\TO}{{  T(\mathcal O)           }}
\nc{\LO}{{  L(\mathcal O)           }}

\nc{\gk}{{      \mathfrak g_\KK             }}
\nc{\bk}{{      \mathfrak b_\KK             }}
\nc{\tk}{{      \mathfrak b_\KK             }}
\nc{\nk}{{      \mathfrak n_\KK             }}
\nc{\go}{{      \mathfrak g_\OO             }}

\nc{\Ll}{{  L_\lambda   }}
\nc{\Lm}{{  L_\mu       }}

\nc{\Gl}{{      {\mathcal G_\lambda}        }}
\nc{\Glb}{{     \barr{\Gl}              }}
\nc{\Glm}{{     \GG_{\lambda+\mu}       }}
\nc{\Ge}{{      \GG_\eta                }}
\nc{\Geb}{{     \barr{\Ge}              }}
\nc{\Gwl}{{     \GG_{W \lambda}         }}
\nc{\Gn}{{      \GG_\nu                 }}
\nc{\Gnb}{      \barr{\GG_\nu           }}
\nc{\Gm}{{      \GG_\mu                 }}
\nc{\Gum}{{     \GG^\mu         }}
\nc{\Gul}{{     \GG^\lambda     }}
\nc{\Gun}{{     \GG^\nu         }}
\nc{\Pl}{{  \PP_\lambda         }}
\nc{\Pwl}{{  {\PP}_{W \lambda}  }}
\nc{\FFl}{{  {\FF}_\lambda      }}

\nc{\Cn}{{  C_\nu       }}
\nc{\Ce}{{  C_\eta      }}

\nc{\Sn }{{     S_\nu           }}
\nc{\Snb}{{  \barr{S_\nu}       }}
\nc{\Sl }{{  S_\lambda          }}
\nc{\Sm }{{  S_\mu              }}

\nc{\GGG}{{  G(\KK)\underset{G(\OO)}\times      \GG                     }}
\nc{\GGlm}{{  (G(\KK)\underset{G(\OO)}\times    \GG)_{(\lambda,\mu)}    }}

\nc{\dgg}{{     \ddot\GG                }}
\nc{\dlm}{{     \ddot\GG_{\la,\mu}      }}

\nc{\db}{{      \hat \Delta     }}
\nc{\dr}{{      \Delta_\R       }}
\nc{\Dp}{{      \Delta^+        }}
\nc{\Dm}{{      \Delta^-        }}
\nc{\da}{{      \Delta_a        }}
\nc{\dap}{{     \Delta_a^+      }}
\nc{\di}{{      \Delta_I        }}
\nc{\cro}{{     \check \rho     }}      
\nc{\dc}{{      \check \delta   }}
\nc{\rac}{{     \check \rho_a   }}


\nc{\Irr}{\operatorname {Irr}}
\nc{\td}{\widetilde d}
\nc{\tg}{\widetilde g}
\nc{\tv}{\widetilde v}
\nc{\tp}{\tilde p}

\nc{\TD}{\widetilde D}
\nc{\TV}{\widetilde V}
\nc{\TU}{\widetilde U}
\nc{\TT}{\widetilde T}
\nc{\TA}{\widetilde A}
\nc{\TB}{\widetilde B}
\nc{\TN}{\widetilde{\mathcal N}}

\nc{\tga}{\widetilde\ga}
\nc{\tde}{\widetilde\de}
\nc{\tphi}{\widetilde\phi}

\nc{\cla}{{\check\lambda}}
\nc{\cmu}{\check\mu}

\nc{\mmm}{{\bf m}}

\nc{\Sym}{\operatorname {Sym}}
\nc{\Spec}{\operatorname {Spec}}
\nc{\Id}{\operatorname {Id}}

\nc{\gl}{\mathfrak{gl}}

\nc{\Perv}{\operatorname {Perv}}
\nc{\Rep}{\operatorname {Rep}}
\nc{\Mlt}{\operatorname {Mlt}}

\nc{\IC}{\operatorname {IC}}

\nc{\thh}{	^{\text{th}}	}                     	
\nc{\subb}{{	\supseteq	}}         %
\nc{\tii}{\widetilde }
\nc{\df}{{ \protect\overset{ \text{def}}= }}		
\nc{\inv}{{ {}^{-1}      }}			


\title{
On quiver varieties
and 
affine Grassmannians of 
\nolinebreak 
type 
\nolinebreak 
$A$
}

\author{       Ivan Mirkovi\'c                 }
\address{\tiny Dept. of Mathematics and Statistics, University
of Massachusetts at Amherst, Amherst MA 01003-4515, USA}
\email{      mirkovic@math.umass.edu   }

\author{       Maxim Vybornov               } 
\email{                vybornov@math.umass.edu }

\begin{abstract} 
We construct Nakajima's quiver varieties of type 
$A$ 
in terms
of affine Grassmannians of type 
$A$. 
This gives a compactification of quiver varieties  
and a decomposition of affine Grassmannians into a 
disjoint union of quiver varieties.
Consequently, 
singularities of
quiver varieties, nilpotent orbits and 
affine Grassmannians are the same in type $A$.
The construction also provides a geometric framework for
skew $(GL(m),GL(n))$ duality and 
identifies the natural basis of
weight spaces in Nakajima's construction
with the
natural basis of multiplicity spaces in tensor products
which arises from affine Grassmannians.
\end{abstract}


\maketitle

\begin{center}
\emph{Dedicated to Igor Frenkel on the occasion of his 50-th birthday}
\end{center}

\se{\bf
Preliminaries
}

\sus{Quiver varieties of type A}
\label{prelimquiver}
We recall Nakajima's construction of simple representations of 
$SL(n)$, cf. \cite{N94, N98}.
Let $I=\{1,\dots,n-1\}$ be the set of vertices
and $H$ be the set of arrows 
of the Dynkin quiver of type $A_{n-1}$.
For an arrow $h\in H$ we denote
by $h'\in I$  and  $h''\in I$ 
its initial and terminal 
vertices.
For a pair $v,d$ in $\Z_{\geq 0}^I$
take $\C$-vector spaces $V_i$ and $D_i$
of dimensions $\dim V_i=v_i$ and $\dim D_i=d_i$, $i\in I$.
Consider the affine space
$$
M(v,d)=\bigoplus_{h\in H}\Hom(V_{h'},V_{h''})\oplus
\bigoplus_{i\in I}\Hom(D_{i},V_{i})\oplus
\bigoplus_{i\in I}\Hom(V_{i},D_{i})
$$
with the natural action of the group 
$G(V)=\prod_{i\in I}GL(V_i)$.
Let $\mmm:M(v,d)\to \mathfrak{g}(V)$ be the corresponding moment map
into the Lie algebra $\mathfrak{g}(V)$. 
Denote $\La(v,d)=\mmm^{-1}(0)$.

Nakajima's quiver variety $\fM(v,d)$ 
is the geometric quotient of $\La^s(v,d)$ by $G(V)$,
where $\La^s(v,d)$ is the set of all stable 
points in 
$\La(v,d)$ (so $\La^s(v,d)/G(V)$ is the set of $\C$-points of
$\fM(v,d)$).
The quiver variety $\fM_0(v,d)=\La(v,d)//G(V)$ is the
invariant theory quotient 
(the spectrum of the $G(V)$-invariant functions).
There is a natural projective map $p:\fM(v,d)\to \fM_0(v,d)$, 
cf. \cite{N98}, and following Maffei \cite{M}, denote 
its image by $\fM_1(v,d)=p(\fM(v,d)) \sub \fM_0(v,d)$.
Finally, let  $\fL(v,d)\df p^{-1}(0)\sub \fM(v,d)$ and 
denote by
$
\HH(\fL(v,d))
$ 
its  top-dimensional Borel-Moore homology. 

\Theo{}
\label{theonakajima}
\cite[\emph{10.ii}]{N98} 
The space $\pl_v \HH(\fL(v,d))$ has the structure of a 
simple $SL(n)$-module with the highest weight $d$
(i.e., $\sum_I\ d_i\om_i$ for the  fundamental weights 
$\om_i$).
The summand
$\HH(\fL(v,d))$ is the weight space for
the weight
$d-Cv$, where $C$ is the Cartan matrix of type $A_{n-1}$.

\sus{From $SL(n)$ to $GL(n)$}
\label{sltogl}
We may consider $\pl_v \HH(\fL(v,d))$ as a representation $W_{\cla}$ 
of $GL(n)$
with highest weight $\cla$,
where $\cla=\cla(d)=(\cla_1, \cla_2,\dots,\cla_n)$ 
is a partition of 
$N=\sum_{j=1}^{n-1}jd_j$ defined as follows:
$\cla_i=\sum_{j=i}^{n}d_j$ (here $d_n=0$). Then
$\HH(\fL(v,d))$ is the weight space $W_{\cla}(a)$,
where $a_i=v_{n-1}+\sum_{j=i}^{n}(d-Cv)_j$ (here $(d-Cv)_n=0$),
cf. \cite[8.3]{N94}.

\sus{Affine Grassmannians of type A} 
We recall the construction of representations 
of $G=GL(m)$
in terms of its affine Grassmannian $\GG_G$, cf. \cite{L81,G95,MV}. 
Let  $V$ be a vector 
space  with a basis $\{e_1,\dots,e_m\}$
and $V((z))=V\otimes_{\C} \C((z))\subb\ L_0=V\otimes_{\C} \C[[z]]$. 
A \emph{lattice} in $V((z))$ 
is an $\C[[z]]$-submodule $L$ of $V((z))$
such that $L\otimes_{\C[[z]]} \C((z))=V((z))$.
The affine Grassmannian $\GG_{G}$ is an ind-scheme 
whose $\C$-points can be described as all lattices in $V((z))$
or as $G((z))/G[[z]]$.
Its connected components $\GG_{(N)}$
are indexed by integers $N\in \Z$, and if $N\ge 0$
then $\GG_{(N)}$ contains
$ 
\GG_N=\{\text{lattices } L \text{ in } V((z)) \text{ such that }
L_0\subseteq L, \dim L/L_0=N \}.
$
To a dominant coweight $\la\in\Z^m$ of $G$, one attaches 
the lattice $L_{\la}=\ \pl_1^m\ \C[[z]]\cd z^{-\la_i}e_i$.
The $G[[z]]$-orbits $\GG_\la$ in $ \GG_{G}$
are parameterized by the dominant coweights $\la$ 
via $\GG_\la=G[[z]]\cd L_\la$.
Finally, we denote by
$L^{<0}G$ the congruence subgroup
of the group ind-scheme
$G[z^{-1}]$ i.e., the kernel of the evaluation
$z^{-1}\mapsto 0$.

The intersection homology of the closure
$\barr\GG_{\la}$ is a realization of the representation $V_{\la}$,
and the convolution of $\IC$-sheaves corresponds to the
tensor products of representations, cf. \cite{G95, MV}.

\sus{Resolution of singularities}
The closure $\barr\GG_{\mu}$ of the orbit
$\GG_{\mu}$ in $\GG_N$ has a natural resolution. 
The $G[[z]]$-orbits in $\GG_{N}$ correspond to $\mu$'s
which may be considered as partitions of $N$ (into at most $m$ parts). 
Any
permutation 
$a=(a_1,\dots,a_n)$
of the partition 
$\cmu$
dual to $\mu$ 
defines a convolution space 
$\widetilde\GG^a_{\mu} =
\GG_{\om_{a_1}}\ast\cddd\ast \GG_{\om_{a_n}}$,
where $\om_{k}$ is the $k$-th fundamental coweight of $G$,
and
a resolution of singularities
$\pi=\pi^a_{\mu}:\widetilde\GG^a_{\mu}\ra\barr\GG_{\mu}$, 
cf. \cite{MV}.

\se{\bf
Nilpotent cones of type A
}

\sus{
$n$-flags \cite{G91,CG}
} 
Let us fix a vector space $D$ 
of dimension $N$. Let $\NN=\NN(D)$ be the nilpotent cone in 
$\End(D)$.
The connected components $\FF^{n,a}$ of the variety of 
$n$-step flags in $D$ 
are parameterized by
all $a\in\Z_{\geq 0}^n$ 
such that 
$N=\sum_{i=1}^{n}a_i$\ :
$$
\FF^{n,a} =\{0=F_0\sub F_1\sub F_2\sub 
\dots\sub F_n=D ~|~\dim F_i-\dim F_{i-1}=a_i\}.
$$
Its cotangent bundle is
$\TN^{n,a}=T^*\FF^{n,a}=\{(u,F)\in\NN\times\FF^a ~|~ u(F_i)
\subseteq F_{i-1} \}$.
Denote by $\mmm_{a}:\TN^{n,a}\ra \NN$ the projection onto 
the first factor.

\sus{
A transverse slice to a nilpotent orbit
}\label{slice}
Let $x$ be a nilpotent operator on $D$,
with Jordan blocks of sizes
$\la=(\la_1\geq\la_2\geq\dots\geq \la_m)$. 
We construct a ``transverse slice'' $T_x$ to the nilpotent orbit 
$\OO_{\la}\sub \NN$ at $x$, 
different from the one considered by Slodowy 
\cite[7.4]{S}.
In some basis
$e_{k,i}$, $1\leq k\leq \la_i$, 
of $D$, one has 
$x :e_{k,i}\mapsto e_{k-1,i}$ (we set $e_{0,i}=0$).
Now
$$
T_x\df\{ x+f,\  f\in\End(D) \ |\ 
f^{l,j}_{k,i} =0, \text { if } k\neq \la_i, 
\text{ and }
f^{l,j}_{\la_i,i} =0, \text { if } l>\la_i \},
$$
where $f^{l,j}_{k,i}:\C e_{l,j}\to\C e_{k,i}$
are the matrix elements of $f$ in our basis.
For a larger orbit $\OO_\mu$,
any permutation $a=(a_1,\dots,a_n)$ of the dual partition $\cmu$,
gives a resolution
$
\TT_{x}^{a}\df\mmm_{a}^{-1}(T_x\cap \barr\OO_{\mu})\subset
\TN^{n,a}
$
of the slice $T_{x,\mu}\df T_x\cap \barr\OO_{\mu}
$
to $\OO_\la$ in $\barr\OO_{\mu}$.

\se{\bf
Main theorem
}
\sus{
From quiver data to  
affine Grassmannian data
}
\label{notation}
We start with $A_{n-1}$ quiver data $v,d\in\Z_{\geq 0}^I$ 
such that $\fM(v,d)$ is nonempty.
Take the $SL(n)$-weights $d$ and $d-Cv$,
and pass to $GL(n)$-weights $\cla$ and  $a$
as in subsection \ref{sltogl}. Now
permute $a$ to a partition
$\cmu=\cmu(a)=(\cmu_1\geq\cmu_2\geq\dots\geq\cmu_n)$ 
of $N=\sum_{j=1}^{n-1}jd_j$. Finally, 
let $\la=(\la_1,\dots,\la_m)$
and $\mu=(\mu_1,\dots,\mu_m)$, 
where $m=\sum_{i=1}^{n-1}d_i$,
be the partitions
of $N$ (i.e., $GL(m)$-coweights) dual to $\cla$ and $\cmu$
respectively. 

\Theo{}
\label{theorem}
Let $N,v,d,a,\la,\mu$ be as above.
Let $L_\la\in\GG_G$ be the lattice corresponding to the 
coweight $\la$,
and let $T_\la$ be its $L^{<0}G$-orbit.
There exist algebraic isomorphisms  
$\phi, \tii\phi, \psi, \tii\psi$ 
such that the following diagram commutes:
\begin{equation}
\begin{CD}
\fM(v,d)  
@>{\tii\phi}>{\simeq}> 
\TT_{x}^{a} 
@>{\tii\psi}>{\simeq}>
\pi^{-1}(T_\la\cap\barr\GG_{\mu}) 
@>{\sub}>> 
\widetilde\GG^{a}_{\mu}  
\\
@V{p}VV  
@V{\mmm_a}VV 
@V{\pi}VV  
@V{\pi}VV 
\\
\fM_1(v,d) 
@>{\phi}>{\simeq}> 
T_{x,\mu} 
@>{\psi}>{\simeq}>
T_\la\cap\barr\GG_{\mu}  
@>{\sub}>> 
\barr\GG_{\mu} 
\end{CD} 
\end{equation}
and $(\psi\circ\phi)(0)=L_\la$.
In particular, 
$\tii\psi\circ\tii\phi$ 
restricts to an isomorphism 
$
\fL(v,d)\simeq\pi^{-1}(L_\la)$.

\sus{} For $d=(d_1,0,\dots,0)$ and $\la=(1,\dots,1)$
the theorem above was proven in (or follows immediately from)
\cite{L81, N94}. 
The isomorphisms $\phi$  (resp. $\tii\phi$) is analogous to the
isomorphism constructed in \cite{N94} 
(resp. isomorphism conjectured in \cite[8.6]{N94} 
and constructed in \cite{M} using a result from \cite{L98}). However, 
our isomorphism $\phi$ is given by an explicit formula
described as follows.
Let us think of a point in $\fM_1(v,d)$ as (closed orbit of) 
a quadruple
$(\{B_i\}_{i\in I},\{\barr{B}_i\}_{i\in I}, 
\{p_i\}_{i\in I}, \{q_i\}_{i\in I})\in\La(v,d)$, 
where $B_i\in \Hom(V_{i},V_{i+1})$, 
$\barr{B}_i\in \Hom(V_{i+1},V_{i})$, 
$p_i\in \Hom(D_{i},V_{i})$, and
$q_i\in \Hom(V_{i},D_{i})$. We decompose the vector space $D$, 
$\dim D=N$, as a direct sum: 
$D=\oplus_{1\leq h\leq j\leq n-1} D_{j}^{h}$, cf. \cite{M},
where $D_{j}^{h}=\C\{e_{h,i}\ | \ \la_i=j \}$, 
$\dim D_{j}^{h}=\dim D_j=d_j$
(notation of \ref{prelimquiver}, 
\ref{slice}, \ref{notation}). Then for any $f\in \End(D)$ 
we consider its blocks
$f^{j',h'}_{j,h}: D_{j'}^{h'}\to D_{j}^{h}$.
By definition, $\phi(B,\barr{B},p,q)=x+f\in \NN$ 
(notation of \ref{slice}), where
\begin{equation}\label{dega}
f^{j',h'}_{j,h}=
\begin{cases}
q_jB_{j-1}\dots B_{h'+1}B_{h'}\barr{B}_{h'}\barr{B}_{h'+1}\dots  
\barr{B}_{j'-1}p_{j'},& \text{ if } h=j, \\
0, & \text{ otherwise. }
\end{cases} 
\end{equation}
In particular, $\phi(0)=x$.

\sus{Compactification of quiver varieties} 
A compactification of $\fM_1(v,d)$ and $\fM(v,d)$
is given by closures of their respective images
under the embeddings
$\fM_1(v,d)\hookrightarrow\barr\GG_{\mu}$ 
and $\fM(v,d)\hookrightarrow\widetilde\GG^{a}_{\mu}$.

\sus{Decomposition} 
The theorem implies a
decomposition of $\barr\GG_{\mu}$ into a  disjoint union
of quiver varieties
\begin{equation}\label{decomp}
\barr\GG_{\mu}=\bigsqcup_{\GG_\la\sub\barr\GG_{\mu}}\ \ 
\bigsqcup_{y\in G\cdot L_\la}
\fM_0(v,d)_y, 
\end{equation}
where $\fM_0(v,d)_y$ is a copy of quiver variety $\fM_0(v,d)$
for every point $y\in G\cdot L_{\la}$, and
$v,d$ are obtained from $\la,\mu$ 
by reversing formulas in subsection \ref{sltogl}.

\sus{
Beilinson-Drinfeld 
Grassmannians
}\label{bd} 
Recall the moment map 
$\mmm:M(v,d)\to \mathfrak{g}(V)$ from 
subsection \ref{prelimquiver}.
Any   $c=(c_1\Id_{V_1},\dots,c_{n-1}\Id_{V_{n-1}})$ 
in the center of the Lie algebra $\mathfrak{g}(V)$
defines  $\La_c(v,d)=\mmm^{-1}(c)$, and then, as in
\ref{prelimquiver},  
the ``deformed'' quiver varieties
$\fM^c(v,d)=\La^s_c(v,d)/G(V)$
and  $\fM^c_{0}(v,d)=\La_c(v,d)//G(V)$.
We expect that in type A our theorem 
and decomposition (\ref{decomp})
extend to a relation between
deformed quiver varieties and the 
Beilinson-Drinfeld 
Grassmannians, cf. \cite{BD}.

For instance, when
$d=(d_1,0,\dots,0)$ 
there is an embedding
$
\fM^c_{0}(v,d)\hookrightarrow \GG^{BD}_{\A^{(n)}}(GL(m))
$
of our quiver variety into the fiber of the
Beilinson-Drinfeld Grassmannian  over the point
$(0,c_1,c_1+c_2,\dots,c_1+\dots+c_{n-1})\in \A^{(n)}$.

The proofs and more details will appear in a forthcoming paper.

Another example of a decomposition 
of an infinite Grassmannian into a disjoint union of quiver varieties
can be found  in \cite{BGK}
(who generalized a result from \cite{W}).
A part of adelic Grassmannian
is  a  union of quiver varieties
$\fM^c(v,d)$ associated to 
{\em affine} 
quivers of type A.

\se{\bf
Geometric construction of skew $(GL(n),GL(m))$ duality
}

\sus{Skew Howe duality} Let $V=\C^m$ and $W=\C^n$ be 
two vector spaces.
Then we have the $GL(V)\times GL(W)$-decomposition
\cite[4.1.1]{H}:
\begin{equation}\label{howew}
\wedge^N(V\otimes W)=\bigoplus_{\la}V_\la\otimes W_{\cla} ,
\end{equation}
where $\la$ varies over all partitions of $N$ 
which fit into the $n\times m$
box, and $V_\la$ and  $W_{\cla}$
are the corresponding
highest weight representation of $GL(m)$ of  $GL(n)$.
This is essentially equivalent to 
natural isomorphisms of vector spaces
\begin{equation}
\label{weightmult1}
\Hom_{GL(m)}
(\wedge^{a_1}V\otimes\dots\otimes \wedge^{a_n}V,
V_{\la})\simeq W_{\cla}(a),
\end{equation}
where $W_{\cla}(a)$ is the weight space corresponding to the weight 
$a=(a_1,\dots,a_n)$.


\sus{} We construct 
a based version of the isomorphism (\ref{weightmult1}), 
i.e., a geometric skew Howe duality.
More precisely, with $N,v,d,a,\la$ as in \ref{notation},
we identify
the right hand side
with
$\HH(\pi\inv (L_\la))$ (notation from Theorem \ref{theorem})
and the left hand side
with
$\HH(\fL(v,d))$ by Theorem\ref{theonakajima}.
The identification of irreducible components
$\Irr\pi\inv (L_\la)=\Irr\fL(v,d)$, 
which follows from Theorem \ref{theorem},
matches the 
natural basis of the space
of intertwiners 
$\Hom_{GL(m)}
(\wedge^{a_1}V\otimes\dots\otimes \wedge^{a_n}V,
V_{\la})$ 
arising from the affine Grassmannian construction 
(i.e., $\Irr\pi\inv (L_\la)$), and the
natural basis of the weight space $W_{\cla}(a)$ in the Nakajima 
construction
(i.e., $\Irr\fL(v,d)$). Altogether:
$$
\Hom_{GL(m)}
(\wedge^{a_1}V\otimes\dots\otimes \wedge^{a_n}V,V_{\la})
\simeq\HH(\pi\inv (L_\la))
\simeq\HH(\fL(v,d))
\simeq W_{\cla}(a).
$$

\sus*{Acknowledgments} 
We are grateful to A. Braverman, D. Gaitsgory, V. Ginzburg, 
M. Finkelberg, A. Maffei, and A. Malkin for useful discussions,
and to MSRI for their 
hospitality and support. The research of M.V. was 
supported by NSF Postdoctoral 
Research Fellowship.

\nopagebreak


\vspace{-.1in}

\begin{thebibliography}{BGK}

{\tiny
\bibitem[BGK]{BGK}
V.~ Baranovsky, V~ Ginzburg, and A. Kuznetsov,
\emph{Wilson's Grassmannian and a noncommutative quadric},
preprint, 2002, math.AG/0203116.

\bibitem[BD]{BD}
A.~ Beilinson and V.~ Drinfeld,
\emph{Quantization of Hitchin's integrable system
and Hecke eigensheaves}, preprint.

\bibitem[CG]{CG}
N.~ Chriss and V.~Ginzburg, 
\emph{Representation theory and complex geometry}, 
Birkh\" auser, Boston, 1997.

\bibitem[G1]{G91}
V.~Ginzburg,  
\emph{Lagrangian construction of the enveloping algebra 
$U({\rm sl}\sb n)$}, 
C. R. Acad. Sci. Paris S\' er. I Math. 312 (1991), no. 12, 907--912.

\bibitem[G2]{G95} 
V.~Ginzburg, 
\emph{Perverse sheaves on a loop group and
Langlands duality}, 
preprint, 1995, alg-geom/9511007.

\bibitem[H]{H}
R. ~Howe, 
\emph{Perspectives on invariant theory: Schur duality, 
multiplicity-free actions and beyond}, 
The Schur lectures (1992) (Tel Aviv), 1--182, Israel Math.
Conf. Proc., 8, Bar-Ilan Univ., Ramat Gan, 1995. 

\bibitem[L1]{L81}
G.~ Lusztig, \emph{Green polynomials and singularities of unipotent 
classes}, 
Adv. in Math. 42 (1981), no. 2, 169--178. 

\bibitem[L2]{L98}
G.~ Lusztig, \emph{On quiver varieties}
Adv. in Math. 136 (1998), 141-182.

\bibitem[M]{M}
A.~ Maffei, \emph{Quiver varieties of type A}, preprint, 2000, 
math.AG/9812142.

\bibitem[MV]{MV}
I.~ Mirkovi\'c, and K.~ Vilonen, 
\emph{Perverse sheaves on affine Grassmannians and Langlands duality}, 
Math. Res. Lett. 7 (2000), no. 1, 13--24. 

\bibitem[N1]{N94}
H.~Nakajima, \emph{Instantons on ALE spaces, quiver varieties, 
and Kac-Moody algebras}, 
Duke Math. J. 76 (1994), no. 2,365--416.

\bibitem[N2]{N98}
H.~Nakajima, \emph{Quiver varieties and Kac-Moody algebras},
Duke Math. J. 91 (1998), no. 3, 515--560. 

\bibitem[S]{S}
P.~ Slodowy,
\emph{Simple singularities and simple algebraic groups},
LNM, 815. Springer, Berlin, 1980.

\bibitem[W]{W}
G. ~ Wilson, \emph{Collisions of Calogero-Moser particles and an
adelic Grassmannian. With an appendix by I.~G.~ Macdonald},
Invent. Math. 133 (1998), no. 1, 1--41. 

}
\end{thebibliography}
\end{document}